\documentclass[12pt]{article}
\usepackage{latexsym,amssymb,amsmath}
\usepackage{amsthm, amstext}
\usepackage{array, amsfonts, mathrsfs}
\usepackage{mathrsfs}
\usepackage{graphicx,color}
\usepackage[cp1251]{inputenc}
\usepackage[T2A]{fontenc}

\newtheorem{Theorem}{Theorem}

\newcommand{\mG}{{\mathscr G}}
\newcommand{\mH}{{\mathscr H}}
\newcommand{\mM}{{\mathscr M}}

\newcommand{\fM}{{\mathfrak M_+}}

\begin{document}

\author{M. I. Belishev\thanks {
St. Petersburg Department of V. A. Steklov Institute of Mathematics of the Russian Academy of Sciences,
27 Fontanka, St. Petersburg, 191023 Russia;
belishev@pdmi.ras.ru.
The work was supported by the grant RFBR 18-01-00269A and by the VW Foundation.},
\
S. A. Simonov\thanks {
St. Petersburg Department of V. A. Steklov Institute of Mathematics of the Russian Academy of Sciences,
27 Fontanka, St. Petersburg, 191023 Russia;
St. Petersburg State University, 7--9 Universitetskaya nab., St. Petersburg, 199034 Russia;
sergey.a.simonov@gmail.com.
The work was supported by the grants RFBR 19-01-00565a, 18-31-00185mol\_a, 17-01-00529a.}}

\title{
On an evolutionary dynamical system of the first order with boundary control}

\date{}

\maketitle

\begin{abstract}
The work is carried out as part of the program to construct a new
functio\-nal (so-called {\it wave}) model of symmetric operators. It
is shown that an abstract evolutionary dynamic system of the first
order (with respect to time) with boundary control, which is determined by
a symmetric operator $L_0:{\mathscr H}\to{\mathscr H}$, is
controllable if and only if $L_0$ has no maximal symmetric
parts in~${\mathscr H}$.
\end{abstract}

\noindent{\bf Keywords:} symmetric operator, wave model, maximal part, invariant subspace.

\noindent{\bf AMS MSC:}\,47B25, 34K30, 47D03

\subsection*{0. Introduction}\label{Sec 1 Introduction}

$\bullet$\,\,\,
We consider an abstract dynamical system, whose evolution in the Hilbert space
${\mathscr H}$
is described by the initial-boundary value problem
\begin{align}
    \mathbb{}\label{Eq 1 dyn BC intr} & iu_t+L^*_0u=0,\qquad t>0 \,;\\
    \label{Eq 2 dyn BC intr} & u|_{t=0}=0\,;\\
    \label{Eq 3 dyn BC intr} & \Gamma_1 u=f\,,\qquad t\geqslant 0\,.
\end{align}
Here
$L_0:{\mathscr H}\to{\mathscr H}$
is a closed densely defined symmetric operator with equal nonzero deficiency indices, which is a part of the {\it Green system} -- a set
$\{{\mathscr H}, {\mathscr B}; L_0, \Gamma_1, \Gamma_2\}$,
in which $\mathscr H$ and $\mathscr B$
are complex separable Hilbert spaces; $\Gamma_{1,2}:\mathscr H\to\mathscr B$ are {\it boundary operators}.

These objects are related by the {\it Green formula}
\begin{equation*}
    (L_0^*u, v)_\mathscr H-(u,L_0^*v)_\mathscr H=(\Gamma_1 u, \Gamma_2
    v)_\mathscr B - (\Gamma_2 u, \Gamma_1 v)_\mathscr B\,;
\end{equation*}
$\mathscr B$-valued function
$f=f(t)$
is called the {\it boundary control}, $u=u^f(t)$ is the solution
({\it
trajectory})
in
${\mathscr
H}$,
correctly defined for an appropriate class of controls
$\mathscr M$.
The sets
${\mathscr U}^T=\{u^f(T)\,|\,f \in
\mathscr M\}$
are called {\it reachable}; ${\mathscr U}={\rm
span}\{{\mathscr U}^T\,|\,\, T>0\}$ is the {\it
total reachable set
}.
The system
(\ref{Eq 1 dyn BC intr})--(\ref{Eq
3 dyn BC intr})
is called controllable, if
$\overline{\mathscr U}={\mathscr
H}$.
\smallskip

\noindent$\bullet$\,\,\,
By
$n_{\pm}[A]={\rm
dim\,\,Ker\,}[A^*\pm i{\bf 1}]$
we denote the deficiency indices of the symmetric operator
$A$.
Let us define the class of operators
$$
\fM\,:=\,\{A\,|\,\,A\subset A^*,\,\,n_+[A]=0\}\,.
$$
We call a subspace
$\mathscr G\!\subset\!\mathscr H$
an invariant subspace of the operator
$A$,
if
${\overline{{\rm Dom\,}A\!\cap\!\mathscr G}\!=\!\mathscr G}$
and
${A[{\rm Dom\,}A\!\cap\!\mathscr G]\!\subset\!\mathscr G}$.
Then the operator
$A_\mathscr G:=A|_{\mathscr G}:\mathscr G\to \mathscr G$
is called the part of
$A$ (in the subspace
$\mathscr G$).
\smallskip

The main result of the present paper is the following. The system
(\ref{Eq 1 dyn BC intr})--(\ref{Eq 3
dyn BC intr})
is controllable, if and only if the operator
$L_0$ does not have parts from the class
$\fM$.
\smallskip

This result is analogous to the result established in \cite{BD_2}
for systems of the second order with lower semibounded operator
$L_0$: for them controllability is equivalent to absence of {\it
self-adjoint} parts. Problems for second order systems can be
reduced to problems for first order systems \cite{Bir Sol}, and in
this sense the result of the present paper is stronger. Besides
that our considerations allow to consider a wide class of Green
systems, while in \cite{BD_2} a canonical Green system is used,
which corresponds to the well-known Vishik decomposition for ${\rm
Dom\,}L_0^*$.
\smallskip

\noindent$\bullet$\,\,\,
The work was carried out as part of the program for development of a new functional
(so-called {\it wave}) model of symmetric operators \cite{JOT,BSim_AA, BSim_FAN,BSim_MatSbor}.
The reason for the study of controllability is an assumption that the wave model describes not
the operator $L_0$ itself, but its {\it wave part} acting in the subspace $\overline{\mathscr U}$.
For a number of well-known examples this assumption is justified: with controllability,
the wave model turns out to be essentially identical to the operator itself \cite{BSim_AA, BSim_MatSbor}.
However, there is a reason to believe that this is not always the case and, moreover,
this is obviously not true if controllability is absent. This work prepares the
study of such effects.

\subsection*{1. Green system}\label{Sec 1 Green system}
$\bullet$\,\,\,
We call a collection
$$\{{\mathscr H},
{\mathscr B}; L_0, \Gamma_1, \Gamma_2\}
$$
a {\it Green system}, if $\mathscr H$ and $\mathscr B$ are complex
separable Hilbert spaces, $L_0$ is a closed densely defined
symmetric operator in $\mathscr H$ with nonzero (possibly
infinite) deficiency indices, $\Gamma_1, \Gamma_2$ are linear
operators from $\mathscr H$ to $\mathscr B$, ${\rm
Dom\,}\Gamma_{1,2}\supset {\rm Dom\,}L^*_0$, for which the
equality (the Green formula)
 \begin{equation}\label{Eq Green Formula}
    (L_0^*u, v)-(u,L_0^*v)=(\Gamma_1 u, \Gamma_2
    v)_\mathscr B - (\Gamma_2 u, \Gamma_1 v)_\mathscr B,\qquad
    u,v\in {\rm Dom\,}L^*_0,
 \end{equation}
holds (where
$(\,\cdot\,,\,\cdot\,)$
and
$(\,\cdot\,,\,\cdot\,)_\mathscr B$
are inner products in
${\mathscr H}$
and
${\mathscr B}$).
The spaces
$\mathscr H$
and
$\mathscr B$
are called the inner and the boundary spaces, respectively, the operator $L_0$
is called the basic operator, $\Gamma_{1,2}$ are the boundary operators.

Let
$A: \mathscr F\to\mathscr G$
be a linear operator,
$\mathscr L$ --
a linear set in
${\mathscr F}$.
By
$A|_{\mathscr L}$
we denote the restriction
$A$
to
${\rm Dom\,}A\cap{\mathscr L}$.

Assume additionally that the Green system defined above
satisfies the following conditions:

\noindent{\bf A.}\,\,\,$\overline{{\rm Ran\,}\Gamma_1}={\mathscr
B}$;

\noindent{\bf B.}\,\,\,$L_0=L^*_0\big|_{{\rm Ker\,}\Gamma_1 \cap
{\rm Ker\,}\Gamma_2}$;

\noindent{\bf C.}\,\,\,$L^*_0\big|_{{\rm Ker\,}\Gamma_1}=: L=L^*$.

\noindent
Below, unless otherwise indicated, they are supposed to hold.
\smallskip

\noindent$\bullet$\,\,\,Note that for every densely defined
symmetric operator $L_0$ with equal deficiency indices there is a
family of Green systems (with properties {\bf A, B, C}), for which
it plays the role of the basic operator. Each of these systems is
determined by some self-adjoint extension of the operator $L_0$:
cf. \cite{MMM}, chapter 7, section 7.1.

\subsection*{2. Free dynamics}\label{Sec 2 Free dynamics}
$\bullet$\,\,\,The Green system has an evolutionary problem in
${\mathscr H}$ associated with it. It has the form
\begin{align}
\label{Eq 1 free dyn}& iv_t+Lv=0,\qquad -\infty<t<\infty\,;\\
\label{Eq 2 free dyn}& v|_{t=T}=y
\end{align}
with finite
$T$
and
$y\in{\mathscr H}$.
Let
$E_\lambda$
be the spectral measure of the (self-adjoint) operator
$L$.
For
$y\in{\rm
Dom\,}L$
this problem has a unique solution
\begin{equation}\label{Eq v^y=e^itLy}
    v=v^{y,T}(t)=v^{y,0}(t-T)=e^{i(t-T)L}y=\int\limits_{\mathbb R}e^{i(t-T)\lambda}dE_\lambda y\,,
\end{equation}
satisfying (\ref{Eq 1 free dyn}), (\ref{Eq 2 free dyn})
in the classical sense \cite{Bir Sol};
at the same time $v^{y,T}(t)\in{\rm Dom\,}L$
for every $t$\, \cite{Bir Sol}.
For $y\not\in{\rm Dom\,}L$
the function $v^{y,T}$,
determined by the right-hand side of (\ref{Eq v^y=e^itLy}),
is understood as a generalized solution of the problem. Denote
\begin{equation*}
    v^y:=v^{y,0}.
\end{equation*}

We call the problem
 (\ref{Eq 1 free dyn})--(\ref{Eq 2 free dyn})
a system with free dynamics. Its evolution is determined by the unitary group
 $e^{i(\,\cdot\,\,-\,T)L}$; $v^y$
is the trajectory of the system, $v^{y,T}(t)$ is the state at the moment $t$.
\smallskip

\noindent$\bullet$\,\,\,
Along with the problem
 (\ref{Eq 1 free dyn})--(\ref{Eq 2 free dyn})
consider the problem
 \begin{align}
\label{Eq 1 free dyn rhs}& iw_t+Lw=g,\qquad -\infty<t<\infty\,;\\
\label{Eq 2 free dyn rhs}& w|_{t=T}=0
 \end{align}
with the right-hand side
 $g\in L_{ 1,\,\rm
loc}\left((-\infty,\infty);{\mathscr H}\right)$.
Its (generalized) solution is
 \begin{equation}\label{Eq free dyn rhs solution}
w=w^{g,T}(t)=\frac{1}{i}\,\int\limits_ T^te^{i(t-s)L}g(s)\,ds=\frac{1}{i}\,\int\limits_ T^t
ds\int\limits_ {\mathbb R}e^{i(t-s)\lambda}\,dE_\lambda\, g(s)\,.
 \end{equation}
For
 $g\in C^1_{\rm
loc}\left((-\infty,\infty);{\mathscr H}\right)$
the solution
 $w^{g,T}$
satisfies
 (\ref{Eq 1 free dyn rhs}), (\ref{Eq 2 free dyn rhs})
in the classical sense. We sketch the proof of this result.

For a continuously differentiable
 $g$
we have
 $$
\frac{d}{ds}\int\limits_ \mathbb R
\frac{1\!-\!e^{i(t-s)\lambda}}{i\lambda}\,dE_\lambda g(s)=\int\limits_ \mathbb
R e^{i(t-s)\lambda}\,dE_\lambda g(s)\!+\!\int\limits_ \mathbb R
\frac{1\!-\!e^{i(t-s)\lambda}}{i\lambda}\,dE_\lambda g^\prime(s).
 $$
Integrating in $s$
from $T$
to $t$,
we get
\begin{align*}
    & -\int\limits_ \mathbb R \frac{1-e^{i(t-T)\lambda}}{i\lambda}\,dE_\lambda
    g(T)=\int\limits_ T^tds\int\limits_ \mathbb R e^{i(t-s)\lambda}\,dE_\lambda
    g(s)\\
    & + \!\int\limits_ T^tds\int\limits_ \mathbb R
    \frac{1\!-\!e^{i(t-s)\lambda}}{i\lambda}\,dE_\lambda
    g^\prime(s)\!=\!iw^{g,T}(t)\!+ \!\int\limits_ T^tds\int\limits_ \mathbb R
    \frac{1-e^{i(t-s)\lambda}}{i\lambda}\,dE_\lambda g^\prime(s),
\end{align*}
which allows to write (\ref{Eq free dyn rhs solution})
in the following form:
\begin{equation}\label{Eq repres w^g}
    w^{g,T}(t)=\int\limits_ \mathbb R
    \frac{1-e^{i(t-T)\lambda}}{\lambda}\,dE_\lambda
    g(T)+\int\limits_ T^tds\int\limits_ \mathbb R
    \frac{1-e^{i(t-s)\lambda}}{\lambda}\,dE_\lambda g^\prime(s)\,.
\end{equation}
This representation is more informative: due to presence of $\lambda$
in the denominators both summands on the right-hand side are elements of ${\rm Dom\,}L$
and, consequently, $w^{g,T}(t)\in{\rm Dom\,}L$
for every $t$.
This makes it possible to check correctness of (\ref{Eq free dyn rhs solution})
by a direct substitution of $w^{g,T}$
in the form (\ref{Eq repres w^g}) into the equation (\ref{Eq 1 free dyn rhs})
and the condition (\ref{Eq 2 free dyn rhs}).
Omitting the corresponding calculations, we note that they are
similar to those given in \cite{BD_2} for a system of the second order.

\subsection*{3. Dynamical system with boundary control}\label{Sec 3 DSBC}
$\bullet$\,\,\,
The Green system has another problem associated with it,
 \begin{align}
\label{Eq 1 dyn BC}& iu_t+L^*_0u=0,\qquad t>0\,;\\
\label{Eq 2 dyn BC}& u|_{t=0}=0\,;\\
\label{Eq 3 dyn BC}& \Gamma_1 u=f\,,\qquad t\geqslant 0\,,
 \end{align}
Here  $f=f(t)$ is a ${\mathscr B}$-valued function of time called
the boundary control. The solution
 $u=u^f(\,\cdot\,)$ is the trajectory of the dynamical system, $u^f(t)$ is its state at the moment $t$.
\smallskip

\noindent$\bullet$\,\,\,
Let us discuss the solvability of the problem
 (\ref{Eq 1 dyn
BC})--(\ref{Eq 3 dyn BC}).
Let
 ${\mathscr N}_{\pm i}={\rm Ker\,}[L_0^*\pm i\bf 1]$
be the deficiency subspaces of the operator
 $L_0$.
Their dimensions
 $n_\pm[L_0]\leqslant\infty$
coincide, because $L_0$ has a self-adjoint extension
 $L$.
The classical solution of the problem is a
 ${\rm Dom\,}L^*_0$-valued function of
 $t$,
 continuous for
 $t\geqslant 0$
and continuously differentiable for $t>0$\,\,\cite{Krein Selim}.
As one can see from (\ref{Eq 3 dyn BC}), a necessary condition of
its solvability is existence of the representation $f=\Gamma_1
\phi(t)$ in terms of a function $\phi$, taking values in ${\rm
Dom\,}L^*_0$. By the von Neumann formula we have ${\rm
Dom\,}L^*_0\ni u=u_0+\phi^++\phi^-$ with  $u_0\in {\rm Dom\,}L_0$
and $\phi^\pm\in {\mathscr N}_{\pm i}$\,; moreover, from the
condition {\bf B} imposed on the Green system above, $\Gamma_1
u_0=0$ holds. Therefore in the ``boundary'' condition (\ref{Eq 3
dyn BC}) one can suppose that
$f=\Gamma_1\left[\phi^+(t)+\phi^-(t)\right]$, which is assumed
from everywhere below.

The substitution $u=p+\phi^++\phi^-$ into (\ref{Eq 1 dyn
BC})--(\ref{Eq 3 dyn BC}) with $L_0^*\phi^\pm=\pm i\phi^\pm$ and
an additional requirement $\phi^{\pm}|_{t=0}=0$ leads to the
problem
\begin{align}
\label{Eq 1 dyn BC for p}& ip_t+L^*_0p=-i\psi,\qquad t>0\,;\\
\label{Eq 2 dyn BC for p}& p|_{t=0}=0\,;\\
\label{Eq 3 dyn BC for p}& \Gamma_1 p=0\,,\qquad t\geqslant 0\,,
\end{align}
where
 $\psi:=\phi^+_t+\phi^+ + \phi^-_t-\phi^-$.
Taking into account (\ref{Eq 3
dyn BC for p})
we have $p\in{\rm Dom\,}L$,
by virtue of which
 $L^*_0p=Lp$
and the problem takes the form
 \begin{align*}
& ip_t+Lp=-i\psi,\qquad t>0\,;\\
& p|_{t=0}=0\,.
 \end{align*}
Using
 (\ref{Eq free dyn rhs solution})
we obtain:
 \begin{equation*}
p(t)= -\int\limits_ 0^t e^{i(t-s)L}\psi(s)\,ds= -\int\limits_ 0^t ds\int\limits_ {\mathbb
R}e^{i(t-s)\lambda}\,dE_\lambda\, \psi(s)\,,
 \end{equation*}
which leads to the representation
\begin{align}\label{Eq dyn BC u^f repres}
    u=u^f(t)=\phi^+(t)+\phi^-(t)-\int\limits_ 0^tds\int\limits_ {\mathbb R}e^{i(t-s)\lambda}\,dE_\lambda\, \psi(s)
\end{align}
with the control $f=\Gamma_1[\phi^++\phi^-]$.

Introduce the class of ``smooth'' controls
 \begin{align*}
{\mathscr
M}:=\left\{f=\Gamma_1[\phi^++\phi^-]\,\big|\,\,\phi^\pm\in
C^2_{\rm loc}\left([0,\infty); {{\mathscr N_{\pm
i}}}\right),\,\,{\rm supp\,}\phi^\pm\subset (0,\infty)\right\};
 \end{align*}
at the same time
 $$
\psi=\phi^+_t+\phi^+ + \phi^-_t-\phi^-\in C^1_{\rm
loc}\left([0,\infty); {\mathscr N_{\pm i}}\right),\quad {\rm
supp\,}\psi\subset(0,\infty)\,.
 $$
Using a trick quite similar to passing from (\ref{Eq free dyn rhs
solution}) to (\ref{Eq repres w^g}) in the system (\ref{Eq 1 free
dyn rhs})--(\ref{Eq 2 free dyn rhs}), one can show that for
$f\in{\mathscr M}$ the function $u^f$ defined by the right-hand
side of (\ref{Eq dyn BC u^f repres}), is a unique classical
solution of the problem (\ref{Eq 1 dyn BC})--(\ref{Eq 3 dyn BC})
with the control $f=\Gamma_1[\phi^++\phi^-]$.

It is easy to see that the right-hand side of (\ref{Eq dyn BC u^f
repres}) makes sense also for less restrictive conditions on $f$
($\phi^\pm$). In this case the function $u^f$ defined by the
right-hand side is considered as a generalized solution of the
problem (\ref{Eq 1 dyn BC})--(\ref{Eq 3 dyn BC}).

\subsection*{4. Controllability}\label{Sec 4 Controllability}
\noindent$\bullet$\,\,\,
For the dynamical system with boundary control the sets of the form
\begin{equation*}
    {\mathscr U}^T:=\{u^f(T)\,|\,\,f \in \mathscr M\}
\end{equation*}
are called reachable (by the moment of time $T$). The set
\begin{equation*}
    {\mathscr U}\,:={\rm span\,}\{{\mathscr U}^T\,|\,\,T>0\}
\end{equation*}
is called the total reachable set. The system (\ref{Eq 1 dyn BC})--(\ref{Eq 3 dyn BC})
is called controllable, if
\begin{equation}\label{Eq def controllability}
    \overline{\mathscr U}\,=\,\mathscr H\,.
\end{equation}
In the case $\overline{\mathscr U}\not=\mathscr H$ the subspace
${\mathscr D}:=\mathscr H\ominus\overline{\mathscr U}$ is called
unreachable.
\smallskip

\noindent${\bullet}$\,\,\, Let us call a subspace ${\mathscr
G}\subset{\mathscr H}$ an invariant subspace of the operator $A$,
if $\overline{{\mathscr G}\cap{\rm Dom\,}A}={\mathscr G}$ and
$A[{\mathscr G}\cap{\rm Dom\,}A]\subset{\mathscr G}$. The
restriction $A|_{{\mathscr G}\cap{\rm Dom\,}A}=:A_{\mathscr G}$,
considered as an operator from ${\mathscr G}$ to ${\mathscr G}$,
is called a part of $A$ in${\mathscr G}$. Note that closedness of
every part of $A$ follows from the closedness of $A$.
\smallskip

\noindent$\ast$\,\,\,We make a remark concerning the terms. The
traditional definition of an invariant subspace for an unbounded
operator does not require the density of $\overline{{\mathscr
G}\cap{\rm Dom\,}A}={\mathscr G}$, cf. \cite{AkhGlaz,MMM}. In our
opinion, this is a drawback, since the operator may have many
``redundant'' invariant subspaces. Indeed, if the operator is such
that ${\rm Dom\,}A\setminus{\rm Dom\,}A^2\not=\varnothing$, then,
by the definition of \cite{AkhGlaz,MMM}, for any set
$y_1,\dots,y_n\subset{\rm Dom\,}A\setminus{\rm Dom\,}A^2$ the
subspace ${\rm span}\{y_1,\dots,y_n, Ay_1,\dots,Ay_n\}$ will be
invariant. Our definition seems to be more natural; at the same
time it allows to get all the desired results.
\smallskip

\noindent$\bullet$\,\,\, Recall that $n_{\pm}[A]={\rm
dim\,\,Ker\,}[A^*\pm i{\bf 1}]$ are the deficiency indices of the
symmetric operator $A$. By the definition adopted in the
Introduction the operator $A$ belongs to the class $\fM$, if
$n_+[A]=0$. We add that the operators of this class are either
self-adjoint (if $n_-[A]=0$) or do not have seld-adjoint
extensions (if  $n_-[A]\not=0$). Such operators are called maximal
\,\,\cite{AkhGlaz,Bir Sol}.

Recall that the free dynamical system and the system with boundary
control are determined by the Green system $\{{\mathscr
H},{\mathscr B}; L_0, \Gamma_1, \Gamma_2\}$, which satisfies the
conditions {\bf A, B, C}.

\begin{Theorem}\label{T 1}
The system {\rm(\ref{Eq 1 dyn BC})--(\ref{Eq 3 dyn BC})} is
controllable, if and only if the operator $L_0$ does not have in
$\mathscr H$ parts from the class $\fM$.
\end{Theorem}

The proof consists of two parts.
\smallskip

{\bf 1.}\,\,{\it Necessity.}\,\,\, Let $L_0$ have a part
 $L_{0\,{\mathscr G}}\in \fM$
in the subspace ${\mathscr G}\subset{\mathscr H}$. We show that
${\mathscr G}\subset{\mathscr D}$, and hence ${\mathscr
D}\not=\{0\}$, which means that the system (\ref{Eq 1 dyn
BC})--(\ref{Eq 3 dyn BC}) is not controllable.
\smallskip

\noindent$\ast$\,\, Choose $y\in{\mathscr G}\cap{\rm
Dom\,}L_0={\rm Dom\,}L_{0\,{\mathscr G}}$ and consider the system
{\it in the subspace } $\mG$:
 \begin{align}
\label{Eq 1 free dyn auxil}& iv_t+L_{0\,{\mathscr G}}v=0,\qquad t<0\,;\\
\label{Eq 2 free dyn auxil}& v|_{t=0}=y\,.
 \end{align}
Since $n_+[L_{0\,{\mathscr G}}]=0$, the half-plane $\mathbb C_+$
consists of regular points of the operator
 $L_{0\,{\mathscr G}}$. Consequently, the problem
(\ref{Eq 1 free dyn auxil})--(\ref{Eq 2 free dyn auxil})
has a classical solution
 $v^y_{\mathscr G}(t)\in{\rm Dom\,}L_{0\,{\mathscr G}}$
for $t\leqslant 0$. This fact follows from the Lumer--Phillips
theorem
 \cite[Th. X.48]{Reed Simon 2}
applied to the accretive operator $A=iL_{0\,{\mathscr G}}$, for
which the half-plane $\{\lambda\in\mathbb C\ |\ {\rm
Re\,}\lambda<0\}$ does not contain points of the spectrum. Note
that such a problem is solvable for times
 $t<0$ and may have no solution for
 $t>0$.

Due to $L_{0\,{\mathscr G}}\subset L$ solutions of the problems
(\ref{Eq 1 free dyn auxil})--(\ref{Eq 2 free dyn auxil}) and
(\ref{Eq 1 free dyn})--(\ref{Eq 2 free dyn}) (with $T=0$) satisfy
the same equation and hence for the same $y$ coincide in $\mH$
for every $t\leqslant 0$. Thus for  $y\in{\rm Dom\,}L_{0 \mathscr G}$
the trajectory $v^y$ of the system (\ref{Eq 1 free dyn})--(\ref{Eq 2 free dyn})
lies in $\mG$ and, moreover, does not leave the linear set ${\rm Dom\,}L_{0 \mathscr G}$.
From the last fact owing to ${\rm Dom\,}L_{0\,{\mathscr G}}\subset{\rm Dom\,}L_0$
and the condition {\bf B} it follows that
 \begin{equation}\label{Eq Gamma v^y=0}
\Gamma_1v^y(t)=\Gamma_2 v^y(t)=0,\qquad t\leqslant 0.
 \end{equation}

\noindent$\ast$\,\,Derive an auxiliary relation. Let $y\in{\rm
Dom\,}L$ and $f\in{\mathscr M}$, so that $v^y$ and $u^f$ are
classical solutions of the problems (\ref{Eq 1 free
dyn})--(\ref{Eq 2 free dyn}) and (\ref{Eq 1 dyn BC})--(\ref{Eq 3
dyn BC}). In the calculation below we simplify the notation by
writing $v$ and $u$ instead of $v^{y,T}$ and $u^f$ and omit the
argument $t$. Integrating by parts using ${\Gamma_1 u=f}$,
$\Gamma_1v=0$ and $iv_t+L_0^*v=iv_t+Lv=0$, we have:
 \begin{align*}
&
0=\int\limits_ 0^T(iu_t+L_0^*u,v)\,dt=\int\limits_ 0^T(iu_t,v)\,dt+\int\limits_ 0^T(L_0^*u,v)\,dt=\langle\text{see}\,(\ref
{Eq 2 dyn BC}),(\ref{Eq 2 free dyn}),(\ref{Eq Green
Formula})\rangle\\
& =
i(u(T),v(T))+\int\limits_ 0^T(u,iv_t)\,dt+\int\limits_ 0^T(u,L_0^*v)\,dt+\int\limits_ 0^T(\Gamma_1u,\Gamma_2v)_{\mathscr B}\,dt
\\
&-\int\limits_ 0^T(\Gamma_2u,\Gamma_1v)_{\mathscr B}\,dt=i(u(T),y)+\int\limits_ 0^T(u,iv_t+L_0^*v)\,dt+\int\limits_ 0^T(f,\Gamma_2v)_{\mathscr B}\,dt\\
& = i(u(T),y)+\int\limits_ 0^T(f,\Gamma_2v)_{\mathscr B}\,dt,
 \end{align*}
whence
 \begin{equation}\label{Eq Auxilliary}
(u^f(T),y)=i\int\limits_ 0^T(f(t),\Gamma_2v^{y,T}(t))_{\mathscr B}\,dt\,=i\int\limits_ 0^T(f(t),\Gamma_2v^y(t-T))_{\mathscr B}\,dt\,.
 \end{equation}
$\ast$\,\,Let again $y\in{\mathscr G}\cap{\rm Dom\,}L_0={\rm
Dom\,}L_{0\,{\mathscr G}}$. From (\ref{Eq Gamma v^y=0}) and
(\ref{Eq Auxilliary}) we have $(u^f(T),y)=0$, from which, owing to
arbitrariness of $f$, we conclude that $y\bot{\mathscr U}^T$.
Since $T>0$ is arbitrary too, we come to $y\bot{\mathscr U}$.

The set of elements $y\in{\mathscr G}\cap{\rm Dom\,}L_0={\rm
Dom\,}L_{0\,{\mathscr G}}$ used above is dense in ${\mathscr G}$.
Therefore ${\mathscr G}\bot{\mathscr U}$, which means that
${\mathscr G}\subset{\mathscr H}\ominus\overline{\mathscr
U}={\mathscr D}\not=\{0\}$.
\medskip

{\bf 2.}\,\,{\it
Sufficiency
.}\,\,
Let us show that in the case
${\mathscr D}\not=\{0\}$
the operator
 $L_0$
has a part
$L_{0\,\mathscr D}\in\fM$.
\smallskip

\noindent$\ast$\,\,
Integrating with respect to time in (\ref{Eq v^y=e^itLy}) we easily get
\begin{equation*}
    w=w^{y,T}(t):=\int\limits_ T^tv^{y,T}(s)\,ds=\int\limits_ \mathbb R\frac{e^{i(t-T)\lambda}-1}{i\lambda}\,dE_\lambda y=w^{y,0}(t-T)\,.
\end{equation*}
Presence of  $\lambda$ in the denominator yields inclusion
$w^{y,T}(t)\in{\rm Dom\,}L$ (so that $\Gamma_1w(t)=0$) for every
$t$ and $y\in{\mathscr H}$. Integrating with respect to time in
(\ref{Eq 1 free dyn})--(\ref{Eq 2 free dyn}), we come to the
system
\begin{align}
\label{Eq 1 dyn w}
& iw_t+Lw=iy,\qquad -\infty<t<\infty\,;\\
\label{Eq 2 dyn w} & w|_{t=T}=0\,.
\end{align}
Denote
\begin{equation*}
    w^y:=w^{y,0}.
\end{equation*}
For $f\in{\mathscr M}$, the corresponding classical solution
$u=u^f$ and arbitrary $y\in{\mathscr H}$, taking into account
$L\subset L_0^*,\,\,\Gamma_1w=0$ and $\Gamma_1u=f$, we have:
\begin{align*}
&
\int\limits_ 0^T(iy,u(t))\,dt=\int\limits_ 0^T(iw_t+Lw,u)\,dt=i\int\limits_ 0^T(w_t,u)\,dt+\int\limits_ 0^T(L_0^*w,u)\,dt
\\
& =
i(w,u)\big|^{t=T}_{t=0}+\int\limits_ 0^T(w,iu_t+L_0^*u)\,dt+\int\limits_ 0^T\left[(\Gamma_1
w,\Gamma_2u)_{\mathscr B}-(\Gamma_2w,\Gamma_1u)_{\mathscr B}\right]\,dt
\\
& = - \int\limits_ 0^T(\Gamma_2w,f)_{\mathscr B}\,dt\,.
\end{align*}
It follows that
\begin{equation}\label{Eq aux 1}
\int\limits_ 0^T\left(y,u^f(t)\right)\,dt=i\int\limits_ 0^T\left(\Gamma_2w^{y,T}(t),f(t)\right)_{\mathscr B}\,dt\,.
\end{equation}

\noindent$\ast$\,\, Let $y\in{\mathscr D}$ and hence the integral
on the left-hand side of (\ref{Eq aux 1}) equals zero for every
control $f$. Pick a sequence of functions $\delta_j\in
C^\infty_0(0,T)$ converging to the Dirac measure supported at the
point $t_0\in(0,T)$. Choose $h\in{\rm Ran\,}\Gamma_1$. Putting
$f=\delta_j(t)h\in\mM$ in (\ref{Eq aux 1}) and letting
$j\to\infty$ we get  $0=(\Gamma_2w^{y,T}(t_0),h)_{\mathscr B}$.
Therefore, owing to denseness of ${\rm Ran\,}\Gamma_1$ in
${\mathscr B}$ (see the condition $\bf A$),
$\Gamma_2w^{y,T}(t)=\Gamma_2w^y(t-T)=0$ holds for every
$t\in[0,T]$. Recalling that $w^y(t)\in{\rm Dom\,}L$ (and hence
$\Gamma_1w^y(t)=0$) we see that $w^y(t)\in {\rm Ker\,}\Gamma_1
\cap {\rm Ker\,}\Gamma_2$ for $t\leqslant0$, i.e.,
\begin{equation}\label{Eq w^y in Dom L0} w^y(t)\in{\rm
Dom\,}L_0\,,\quad t\leqslant0\,.
\end{equation}

\noindent$\ast$\,\,
Let us derive another auxiliary relation. For $f\!\in\!\mM$ and ${y\!\in\!\mH}$,
taking into account
$$
iw_t+Lw=iw_t+L_0^*w=iy,\quad\Gamma_1w=0,\quad u|_{t=0}=0\,,
$$
for $T>0$ and $t\leqslant0$, we have:
\begin{align*}
& 0=\int\limits_ 0^T(iu_t(s)+L_0^*u(s),w(s+t))\,ds=(i u(T),w(T+t))
\\
&
+\int\limits_ 0^T\left(u(s),\left[i\frac{d}{dt}+L_0^*\right]w(s+t)\right)\,ds
\\
&+\int\limits_ 0^T\left[(\Gamma_1
u(s),\Gamma_2w(s+t))_{\mathscr B}-(\Gamma_2u(s),\Gamma_1w(s+t))_{\mathscr B}\right]\,ds
\\
& =(i
u(T),w(T+t))+\int\limits_ 0^T(u(s),iy)\,ds+\int\limits_ 0^T(f(s),\Gamma_2w(s+t))_{\mathscr
B}\,ds,
\end{align*}
whence
\begin{equation}\label{Eq aux 2}
(u^f(T),w^y(t))=\int\limits_ 0^T(u^f(s),y)\,ds-i\int\limits_ 0^T(f(s),\Gamma_2w^y(s+t-T))_{\mathscr
B}\,ds\,.
\end{equation}

\noindent$\ast$\,\, Owing to (\ref{Eq aux 2}), from $y\in{\mathscr
D}$, $T>0$ and  $t\leqslant0$ it follows that $w^y(s+t-T)\in{\rm
Ker\,}\Gamma_2$ for every $s\in[0,T]$ and $(u^f(T),w^y(t))=0$.
Therefore $w^y(t)\bot{\mathscr U}$, i.e., the trajectory $w^y(t)$
is contained in ${\mathscr D}$ for $t\leqslant0$. Moreover,
according to  \eqref{Eq w^y in Dom L0}, the set of states
$$
{\mathscr S}:=\left\{w^y(t)\,|\,\,y\in{\mathscr
D},\,t\leqslant0\right\}
$$
is contained in ${\mathscr D}\cap{\rm Dom\,}L_0$.
\smallskip

The set ${\mathscr S}$ is dense in ${\mathscr D}$. Indeed, if
$d\in{\mathscr D}$ and $d\bot{\mathscr S}$, then
$(d,w^y(t))|_{t\leqslant0}=0$ and hence
$(d,w^y_t(t))|_{t\leqslant0}=0$. From the last equality, according
to \eqref{Eq 1 dyn w}--\eqref{Eq 2 dyn w} (with $T=0$), we have
\begin{equation*}
    0=(d,iw^y_t(0))=(d,-Lw^y(0)+iy)=-i(d,y).
\end{equation*}
From arbitrariness of $y\in{\mathscr D}$
we have $d=0$.
\smallskip

Denseness of $\mathscr S$ implies that
\begin{equation}\label{Eq density S}
    \overline{{\mathscr D}\cap{\rm Dom\,}L_0}\supset\overline{\mathscr S}=\mathscr D\,.
\end{equation}
\smallskip

\noindent$\ast$\,\,As shown above, for $y\in\mathscr D$ the
trajectory $w^y(t)$ is contained in $\mathscr D$ for
$t\leqslant0$. Therefore also the trajectory $v^y(t)=w^y_t(t)$ is
contained in $\mathscr D$. If at the same time $y\in{\mathscr
D}\cap{\rm Dom\,}L_0$, then the trajectory $v^y(t)$ is classical
and, according to (\ref{Eq 1 free dyn}), we have
$L_0y=Ly=Lv^y(0)=-iv^y_t(0)\in{\mathscr D}$, i.e.,
\begin{equation}\label{Eq invar}
L_0[{\mathscr D}\cap{\rm Dom\,}L_0]\subset{\mathscr D}\,.
\end{equation}
Comparing (\ref{Eq density S}) and (\ref{Eq invar}) we conclude
that the operator $L_0$ has a (closed symmetric) part
$L_{0{\mathscr D}}=L_0|_{{\mathscr D}\cap{\rm Dom\,}L_0}$ in
${\mathscr D}$.

\smallskip

\noindent$\ast$\,\,Let us show that $L_{0 {\mathscr
D}}\!\in\!\fM$.  To do this, verify the equality ${n_+[L_{0
{\mathscr D}}]\!=\!0}$.

Let $y\in{\mathscr D}\cap{\rm Dom\,}L_0={\rm Dom\,}L_{0{\mathscr
D}}$;  in this case, as mentioned above, for $t\leqslant0$ the
inclusion $v^y(t)\in{\mathscr D}\cap{\rm Dom\,}L_0={\rm Dom\,}L_{0
{\mathscr D}}$ holds. For an element $z\in{\rm Ker\,}[L_{0
{\mathscr D}}^*+i{\bf 1}]=\mathscr N_i[L_{0 {\mathscr D}}]$ of the
defect subspace of $L_{0{\mathscr D}}$ we have:
\begin{align*}
&0=\left([L_{0 {\mathscr D}}^*+i{\bf1}]z,v^y(t)\right)=\left(z,L_{0{\mathscr D}}v^y(t)\right)-(z,iv^y(t))=\\
&=\left(z,Lv^y(t)\right)-(z,iv^y(t))=\langle\,\text{see}\,\,(\ref{Eq 1 free dyn})\,\rangle=\left(z,-iv^y_t(t)\right)-(z,iv^y(t))\,,
\end{align*}
whence
$$
\left(z,v^y(t)\right)_t+\left(z,v^y(t)\right)=0\,,\qquad t\leqslant0.
$$
Using (\ref{Eq 2 free dyn}) we get $(z,v^y(t))=(z,y)e^{-t}$, or
$$
(z,e^{itL}y)=(z,y)\,e^{-t}\,,\qquad t\leqslant0\,.
$$
Since $|(z,e^{itL}y)|\leqslant \|z\|\|y\|$ and the right-hand
side is not bounded, the last equality can hold, only if
$(z,y)=0$. Due to arbitrariness of the choice of $y$ from the set
${\mathscr D}\cap{\rm Dom\,}L_0$, which is dense in ${\mathscr
D}$, we conclude that $z=0$. Therefore $\mathscr N_i[L_{0
{\mathscr D}}]=\{0\}$ and $n_+[L_{0{\mathscr D}}]=0$.
\smallskip

We have shown that the operator  $L_0$ has the part $L_{0\mathscr
D}\in\fM$ in ${\mathscr D}$. Theorem \ref{T 1} is proved.
\smallskip

\noindent$\bullet$\,\,\,
Evolution of the system with boundary control at negative times
is described by the problem
\begin{align}
\label{Eq 1 dyn BC-}& iu_t+L^*_0u=0,\qquad t<0\,;\\
\label{Eq 2 dyn BC-}& u|_{t=0}=0\,;\\
\label{Eq 3 dyn BC-}& \Gamma_1 u=f\,,\qquad t\leqslant 0\,.
\end{align}
Its reachable sets are
\begin{equation*}
    {\mathscr U}^T:=\{u^f(T)\,|\,\,f(-t) \in \mathscr M\}\,, \qquad {\mathscr U}_-\,:={\rm span\,}\{{\mathscr
    U}^T\,|\,\,T<0\}\,.
\end{equation*}
The system (\ref{Eq 1 dyn BC-})--(\ref{Eq 3 dyn BC-}) is controllable (for negative times), if
\begin{equation}\label{Eq def controllability}
    {\mathscr D}_-:=\mathscr H\ominus\overline{\mathscr U}_-=\{0\}\,.
\end{equation}
A proof analogous to the one given above for Theorem \ref{T 1}
leads to a similar result: the system (\ref{Eq 1 dyn
BC-})--(\ref{Eq 3 dyn BC-}) is controllable, if and only if the
operator $L_0$ does not have parts $L_{0\mathscr G}$ in $\mathscr
H$ with the deficiency index $n_-[L_{0 \mathscr G}]=0$.

Let $\mathscr D_+$ and $\mathscr D_-$ be unreachable  subspaces of
the systems (\ref{Eq 1 dyn BC})--(\ref{Eq 3 dyn BC}) and (\ref{Eq
1 dyn BC-})--(\ref{Eq 3 dyn BC-}), respectively. In the
(hypothetical) case $\mathscr D_+\cap\mathscr D_-\not=\{0\}$ the
operator $L_0$ has a self-adjoint part in $\mathscr
D_+\cap\mathscr D_-$.
\smallskip

Recall that symmetric operators that have at least one deficiency
index equal to zero are called maximal. Summarizing our
considerations, we conclude that the system with boundary control
is controllable at {\it all} times (i.e., $\mathscr D_+=\mathscr
D_-=\{0\}$), if and only if the operator $L_0$ has no maximal
parts.
\smallskip

\noindent$\bullet$\,\,\,
As an illustration, consider the Dirac system on the half-line.
\smallskip

Elements of the corresponding Green system are:
\begin{equation*}
\begin{array}{c}
    \mathscr H=L_2(\mathbb R_+;\mathbb C^2);\quad\mathscr B=\mathbb C;\quad L_0=J\frac{d}{dx},
    \\
    {\rm Dom\,}L_0=\{y\in\mathscr H\ |\ y\in AC_{\rm loc}(\mathbb R_+),\ y'\in\mathscr H,\ y(0)=0\};
    \\
    \Gamma_1y=y^1(0),\quad \Gamma_2y=y^2(0),
    \end{array}
\end{equation*}
where
$y=
    \begin{pmatrix}
      y^1 \\
      y^2 \\
    \end{pmatrix}$,
$J=\begin{pmatrix}
   0 & 1 \\
   -1 & 0 \\
\end{pmatrix}$.
At the same time
\begin{equation*}
    L_0^*=J\frac d{dx},\quad {\rm Dom\,}L_0^*=\{y\in\mathscr H\ |\ y\in AC_{\rm loc}(\mathbb R_+),\ y'\in\mathscr H\},
\end{equation*}
and the Green formula is
\eqref{Eq Green Formula}.
\smallskip

The system with  boundary control has the form
\begin{align*}
&iu_t+L^*_0u=0,\qquad x\in\mathbb R_+,\  t>0\,;\\
&u|_{t=0}=0,\qquad x\in\mathbb R_+\,;\\
&\Gamma_1 u=u^1(0,t)=f(t)\,,\qquad t\geqslant 0\,;
\end{align*}
controls are taken from the class $\mathscr M=\{f\in C_{\rm loc}^{\infty}([0,\infty);\mathbb C)\ |\ {\rm supp\,}f\subset(0,\infty)\}$. Its trajectory can be found explicitly: redefining $f|_{t<0}\equiv0$, we have
\begin{equation*}
  u=u^f(x,t)=
  \begin{pmatrix}
    u^1(x,t) \\
    u^2(x,t) \\
  \end{pmatrix}
  =f(t-x)
  \begin{pmatrix}
    1 \\
    i \\
  \end{pmatrix}
  ,\quad x\in\mathbb R_+,\ t\geqslant0.
\end{equation*}
Further, we easily find that
\begin{equation*}
  \overline{\mathscr U^T}=\left\{\varphi
  \begin{pmatrix}
    1 \\
    i \\
  \end{pmatrix}
  \ \Big|\ \varphi\in L_2(\Omega^T;\mathbb C)\right\},
  \quad
  \overline{\mathscr U}=\left\{\varphi
  \begin{pmatrix}
    1 \\
    i \\
  \end{pmatrix}
  \ \Big|\ \varphi\in L_2(\mathbb R_+;\mathbb C)\right\},
\end{equation*}
where $\Omega^T=(0,T)\subset\mathbb R_+$. As a consequence, we obtain that
\begin{equation*}
  \mathscr D=\mathscr H\ominus\overline{\mathscr U}=
  \left\{\psi
  \begin{pmatrix}
    1 \\
    -i \\
  \end{pmatrix}
  \ \Big|\ \psi\in L_2(\mathbb R_+;\mathbb C)\right\}.
\end{equation*}
According to our results, the operator $L_0$ has in $\mathscr D$ the maximal part
\begin{align*}
  L_{0\mathscr D}&=J\frac d{dx},
  \\
  {\rm Dom\,}L_{0\mathscr D}&\!=\!\left\{\psi
  \begin{pmatrix}
    1 \\
    -i \\
  \end{pmatrix}
  \ \Big|\ \psi\in AC_{\rm loc}(\mathbb R_+),\ \psi,\psi'\!\in\! L_2(\mathbb R_+;\mathbb C),\ \psi(0)\!=\!0\right\}.
\end{align*}
It is unitarily equivalent to the scalar operator $i\frac d{dx}$
in $L_2(\mathbb R_+;\mathbb C)$ with the domain
$$
\{f\in H^1(\mathbb R_+;\mathbb C)\ |\ f(0)=0\}.
$$
As a consequence, we have $n_+[L_{0,\mathscr D}]=0$,
$n_-[L_{0,\mathscr D}]=1$ (cf. \cite{Bir Sol}).


\begin{thebibliography}{9}

\bibitem{AkhGlaz}
N. I. Akhiezer, I. M. Glazman.
Theory of Linear Operators in Hilbert Space.
{\em Dover Publications}, 1993.

\bibitem{JOT}
M. I. Belishev.
\newblock{A unitary invariant of a semi-bounded operator in reconstruction of manifolds.}
\newblock{\em J. Operator Theory}, 69(2), 2013, 299--326.

\bibitem{BD_2}
M. I. Belishev, M. N. Demchenko.
Dynamical system with boundary control associated with a symmetric semibounded operator.
{\em Journal of Mathematical Sciences},
194(1), 2013, 8--20.

\bibitem{BSim_AA}
M. I. Belishev, S. A. Simonov.
\newblock {Wave model of the Sturm-Liouville operator on the half-line.}
\newblock{\em St. Petersburg Math. J.},
29(2), 2018, 227--248;
arXiv: 1703.00176v1.

\bibitem{BSim_FAN}
M. I. Belishev, S. A. Simonov.
The wave model of the metric space with measure and its application.
{\em To appear in Sbornik: Mathematics}.

\bibitem{BSim_MatSbor}
M. I. Belishev, S. A. Simonov.
A wave model of metric spaces.
{\em Funct. Anal. Appl.},
53(2), 2019, 79--85.

\bibitem{Bir Sol}
M. S. Birman, M. Z. Solomyak.
Spectral Theory of Self-Adjoint Operators in Hilbert Space.
\newblock{\em D. Reidel Publishing Comp.}, 1987.

\bibitem{MMM}
V. O. Derkach, M. M. Malamud.
Extension theory of symmetric operators and boundary value problems.
Proceedings of Institute of Mathematics of NAS of
Ukraine, vol. 104. {\em Institute of Mathematics of NAS of Ukraine}, 2017.

\bibitem{Krein Selim}
S. G. Krein.
Linear Differential Equations in Banach Space.
{\em AMS}, 1972.

\bibitem{Reed Simon 2}
M. Reed, B. Simon.
Methods of Modern Mathematical Physics. II: Fourier Analysis, Self-adjointness.
{\em Academic Press}, 1975.

\end{thebibliography}
\end{document}